\documentclass[12pt]{article}

\usepackage[utf8]{inputenc}
\usepackage[T1]{fontenc}
\usepackage[margin=1in]{geometry}
\usepackage{amsmath, amsfonts, amssymb, amsthm}
\usepackage[none]{hyphenat}
\usepackage{graphicx}
\usepackage{float}
\usepackage{times}
\usepackage{lipsum}
\usepackage{enumitem}
\usepackage{bbm}

\newcommand{\R}{\mathbb R}
\newcommand{\N}{\mathbb N}

\newtheorem{Th}{Theorem}[section]
\newtheorem{df}{Definition}[section]

\newtheorem{lm}{Lemma}[section]

\def\XXint#1#2#3{{\setbox0=\hbox{$#1{#2#3}{\int}$ }
\vcenter{\hbox{$#2#3$ }}\kern-.6\wd0}}

\numberwithin{equation}{section}

\newcommand{\norme}[1]{\left\Vert #1\right\Vert}

\def\XXint#1#2#3{{\setbox0=\hbox{$#1{#2#3}{\int}$ }
\vcenter{\hbox{$#2#3$ }}\kern-.6\wd0}}

\newcommand\restr[2]{{% we make the whole thing an ordinary symbol
  \left.\kern-\nulldelimiterspace % automatically resize the bar with \right
  #1 % the function
  \vphantom{\big|} % pretend it's a little taller at normal size
  \right|_{#2} % this is the delimiter
  }}

\begin{document}

\title{Mountain pass solutions to equations with subcritical Musielak-Orlicz-Sobolev growth}

\author{Allami Benyaiche and Ismail Khlifi\\
\small {\it allami.benyaiche@uit.ac.ma and is.khlifi@gmail.com}\\
\small{Ibn Tofail University, Department of Mathematics, B.P: 133, Kenitra-Morocco.}}
\date{}
\maketitle

\noindent \textbf{Abstract.} In this paper, we prove the existence of solutions to quasilinear elliptic equations on a bounded domain of $\R^N$ under subcritical Musielak-Orlicz-Sobolev growth. Our proofs rely essentially on Mountain Pass Theorem with corresponding variational techniques. Furthermore, we establish the sharpness of our central assumptions. As far as we know, our approach is new, even for the Orlicz case.

\section{Introduction}

Since the 1970s, the theory of critical points has undergone rapid development in the branch of the calculus of variations because it has wide applications in other mathematical fields. In particular, it is used to prove the existence of solutions to partial differential equations and dynamical systems.\\
\\
The simplest quasilinear elliptic equation is expressible by the p-Laplacian operator
\begin{equation}\label{P p-laplacian}
\left \{
   \begin{array}{rlll}
    \displaystyle -\text{div}\left(|\nabla u|^{p-2} \nabla u\right) & = & f(x,u)  & \text{in} \;  \Omega\\
    u & = & 0 & \text{on} \;  \partial\Omega,
   \end{array}
   \right.
\end{equation}
where $1 < p < \infty$ and $\Omega$ is a bounded domain of $\R^N$. Such problem has been extensively studied in the literature, see for example $\cite{ref1,ref3,ref4}$, because there are some physical phenomena that such kinds of equations can model. One of the approaches most used to prove the existence of solutions to this problem is based on the mountain pass theorem. For the history, A. Ambrosetti and P. Rabinowitz treated the semilinear case ($p=2$) $\cite{ref1}$. next, G. Dinca, P. Jebelean, and J. Mawhin $\cite{ref3}$ used it to solve the problem for general $p > 1$ under the following conditions: $f: \Omega \times \R \to  \R$ is a Carathéodory function that satisfies 
\begin{itemize}
\item The p-subcritical growth condition: 
$$|f(x,t)| \leq C(1 + t^{q-1}) \quad \text{for all} \; t \in \R \; \text{and} \; x \in \Omega,$$
where $p < q <p^* = \frac{Np}{N-p}$ if $p < N$ or $p^* = \infty $ if $p \geq N$. 
\item The p-superlinear at $0$ condition: 
$$\; \limsup_{t \to 0} \frac{f(x,t)}{|t|^{p-2}t} < \lambda_1 \quad \text{uniformly with} \; x \in \Omega,$$
where $\lambda_1$ is the first eigenvalue of p-Laplacian operator with homogeneous Dirichlet boundary on $W^{1,p}_0(\Omega)$.
\item The p-Ambrosetti-Rabinowitz type condition: there exist $\theta > p$ and $t_0 > 0$ such that
$$ 0 < \theta F(x,t) \leq tf(x,t) \quad \text{for all} \; t \in \R \; \text{such that} \; |t| \geq t_0, \; x \in \Omega,$$
where $F(x,t) = \int_0^t f(x,s)  \, \mathrm{d}s$.
\end{itemize}
The restriction to p-subcritical ($p < q <p^*$) ensures that the embedding $W^{1,p}_0(\Omega) \hookrightarrow L^q(\Omega)$ is compact. The Ambrosetti–Rabinowitz condition (AR) is important to ensure that the energy functional associated with the problem has the mountain pass geometry and verify the Palais–Smale condition (PS). Over the years, many researchers trying to weaken the above condition in different situations.\\
\\
Later on, the problem $(\ref{P p-laplacian})$ has been extended to variable exponent case:
\begin{equation}\label{P p(x)-laplacian}
\left \{
   \begin{array}{rlll}
    \displaystyle -\text{div}\left(|\nabla u|^{p(x)-2} \nabla u\right) & = & f(x,u) & \text{in} \;  \Omega\\
  \hspace{3.5cm}  u & = & 0  & \text{on} \;  \partial\Omega
   \end{array}
   \right.
\end{equation}   
where $x \to p(x)$ is a continuous function on $\overline{\Omega}$ such that $1 < p^- := \inf_{x \in \Omega} p(x) \leq p(x) \leq p^+ := \sup_{x \in \Omega}p(x) < \infty$ and $f: \Omega \times \R \to  \R$ is a Carathéodory function satisfies 
\begin{itemize}
\item The p(x)-subcritical growth condition: 
$$|f(x,t)| \leq C(1 + t^{q(x) - 1}) \quad \text{for all} \; t \in \R \; \text{and} \; x \in \Omega,$$
where $p^+ < q^- <  q(x) <p^*(x)= \frac{Np(x)}{N-p(x)}$ if $p(x) < N$ and $p^*(x) = \infty $ if $p(x) \geq N$.
\item The p(x)-superlinear at $0$ condition: 
$$\lim_{t \to 0} \frac{f(x,t)}{|t|^{p^+ - 1}} = 0 \quad \text{uniformly with} \; x \in \Omega.$$ 
\item The p(x)-Ambrosetti-Rabinowitz condition: there exist $\theta > p^+$ and $t_0 > 0$ such that
$$0 < \theta F(x,t) \leq tf(x,t) \quad \text{for all} \; t \in \R \; \text{such that} \; |t| \geq t_0, \; x \in \Omega,$$
where $F(x,t) = \int_0^t f(x,s)  \, \mathrm{d}s$.\\
\end{itemize}
The existence of a nontrivial solution to problem $(\ref{P p(x)-laplacian})$ has been studied in $\cite{ref11}$. For more details in this direction, see $\cite{ref5,ref11}$ and the references therein.\\
\\
Although the operator p(x)-Laplacian is a natural generalization of the p-Laplacian (p(x) = p constant), the variational problem $(\ref{P p(x)-laplacian})$ requires a more complicated analysis which must be performed to study the existence of nontrivial solutions. However, the generalizations of the conditions used in problem $(\ref{P p-laplacian})$ to the exponent variable are not difficult. Yet, it is not always the case as in the Orlicz situation. So, it is a good motivation to study this problem in a general framework as the Musielak-Orlicz-Sobolev case. For this, we consider the following problem
\begin{equation} \label{G()-laplacian}
\left \{
   \begin{array}{rlll}
    \displaystyle -\text{div}\left(\frac{g(x,|\nabla u|)}{|\nabla u|}\nabla u\right) & = & f(x,u) & \text{in} \;  \Omega\\
    u & = & 0 & \text{on} \; \partial\Omega,
   \end{array}
   \right.
\end{equation}
where $g(\cdot,t)$ is the density of a generalized $\Phi$-function $G(\cdot)$ (see section $2$) which satisfy the condition:\\
$(SC)$ There exist two constants $g^0 \geq g_0 > 1$ such that,
$$ g_0 - 1 \leq \displaystyle \frac{t g^\prime(x,t)}{g(x,t)} \leq g^0-1.$$
In the Musielak-Orlicz case, the major difficulty is to give a condition that replaces the subcritical growth in the p-Laplacian and p(x)-Laplacian cases. To overcome this problem, we were inspired by the compact embedding theorem given by P. Harjulehto and P. H\"ast\"o in their fundamental monograph $\cite{ref6}$. So, the main novelty of our work is to introduce a new subcritical assumption: for $g^0 < \min\{g_0^*, N\}$, the G($\cdot$)-subcritical growth condition:
$$|f(x,t)| \leq C(1+\psi(x,t)),$$
with $\psi$ is density of the generalized $\Phi$-function $\Psi(\cdot)$ which satisfies $\Psi^{-1}(x,t) = t^{-\alpha} G^{-1}(x,t)$ with $\alpha \in (\frac{1}{g_0} - \frac{1}{g^0} , \frac{1}{N})$. Such condition coincides with the p(x)-subcritical growth condition for $G(x,t) = t^{p(x)}$. Indeed, the first thing, we have $g_0 = p^-$, $g^0 = p^+$ and $G^{-1}(x,t) = t^{\frac{1}{p(x)}}$. Hence
$$ \Psi^{-1}(x,t) = t^{-\alpha +  \frac{1}{p(x)}} = t^{ \frac{1 - \alpha p(x)}{p(x)}}.$$
So,
$$\Psi(x,t) = t^{\frac{p(x)}{1 - \alpha p(x)}} = t^{q(x)},$$
where $q(x):= \displaystyle  \frac{p(x)}{1 - \alpha p(x)}$. As $\alpha < \displaystyle \frac{1}{N}$, then   
$$ q(x) < \frac{p(x)}{1 - \frac{p(x)}{N}} = \frac{Np(x)}{N- p(x)} = p^*(x). $$
But, we still have to $p^+ < q^-$. For this, we determine the (SC) of $\Psi(\cdot)$. From the condition $(SC)$ of $G(\cdot)$, we have
$$ g_0 \leq \displaystyle \frac{t g(x,t)}{G(x,t)} \leq g^0.$$
By Proposition $2.3.7$ in $\cite{ref6}$ we get
$$ \frac{1}{g^0} \leq \displaystyle \frac{t (G^{-1})^\prime(x,t)}{G^{-1}(x,t)} \leq \frac{1}{g_0}.$$
As 
$$\frac{t (\Psi^{-1})^\prime(x,t)}{\Psi(x,t)} = -\alpha + \frac{t (G^{-1})^\prime(x,t)}{G^{-1}(x,t)},$$
then
$$- \alpha + \frac{1}{g^0} \leq \frac{t (\Psi^{-1})^\prime(x,t)}{\Psi^{-1}(x,t)} \leq -\alpha + \frac{1}{g_0}.$$
Thus 
$$\frac{1}{-\alpha + \frac{1}{g_0}} \leq \frac{t \psi(x,t)}{\Psi(x,t)} \leq \frac{1}{-\alpha + \frac{1}{g^0}}.$$
So, an assumption of type $g^0 < \min\{\frac{1}{-\alpha + \frac{1}{g_0}} , N\}$ is needed. Hence, we restrict that $\alpha \in (\frac{1}{g_0} - \frac{1}{g^0} , \frac{1}{N})$. But, for this last interval to be well defined, we suppose that $g^0 < \min\{g_0^*, N\}$ which implies that $\frac{1}{g_0} - \frac{1}{g^0} < \frac{1}{N}$. For the other two conditions, we can introduce them by the same way as in the variable exponent case (see section $3$). With these assumptions, we are able to establish the existence of nontrivial solutions to problem $(\ref{G()-laplacian})$.

\section{Preliminary results}

We recall some definitions relating to Musielak-Orlicz spaces. A major synthesis of the functional analysis in these spaces is given in the monographs of Musielak $\cite{ref9}$ and Harjulehto, H\"ast\"o $\cite{ref6}$.\\ 
\\
Throughout this work, let $\Omega$ be a bounded domain, $C$ be a generic constant whose value may change between appearances and $g: \Omega \times [0,\infty) \to [0,\infty)$ such that for each $t \in [0,\infty)$, the function $g(\cdot,t)$ is measurable and for a.e. $x \in \Omega$, $g(x, \cdot)$ is a $C^1([0,\infty))$ which satisfies the condition:
\begin{enumerate}
\item[•] $(SC)$ There exist two constants  $g^0 \geq g_0 > 1$ such that,
$$ g_0 - 1 \leq \displaystyle \frac{t g^\prime(x,t)}{g(x,t)} \leq g^0 - 1.$$
\item[•] $(A_0)$ There exists a constant  $C > 1$ such that,
$$ \displaystyle \frac{1}{C} \leq g(x,1) \leq C \, \; \text{a.e.} \; \, x \in \Omega.$$ 
\item[•] $(A_{1})$ If there exists a positive constant $C$ such that, for every Ball $B_R$ with $R < 1$ and $ x,y \in B_R \cap \Omega$, we have
$$ G(x,t) \leq CG(y,t) \; \; \; \text{when} \; \; G^-_{B}(t) \in \left[1 , \frac{1}{R^n}\right].$$
\end{enumerate}
For each $x \in \Omega$ and $t \geq 0$ we define $G(\cdot)$ by
$$ G(x,t) =  \displaystyle \int_{0}^t g(x,s) \, \mathrm{d}s.$$
The function $G(\cdot)$ is a generalized $\Phi$-function, we notice $G(\cdot) \in \Phi(\Omega)$.\\

\noindent \textbf{Properties of generalized $\Phi(\cdot)$-functions.} The first lemma recall some useful inequalities that can be deduced easily from the condition $(SC)$ (see $\cite{ref2,ref8}$).

\begin{lm}{\label{Prs of G}}
Let $G(\cdot) \in \Phi(\Omega)$. Then we have the following properties
\begin{enumerate}[label=(\alph*)]
\item $t \to g(x,t)$ is nondecreasing and
\begin{equation}
ag(x,b) \leq ag(x,a) + g(x,b)b.
\end{equation}
\item $G(\cdot)$ satisfies also the condition $(SC)$
$$ g_0 \leq \displaystyle \frac{t g(x,t)}{G(x,t)} \leq g^0 .$$
\item 
\begin{equation}
\sigma^{g_0} G(x,t) \leq G(x,\sigma t) \leq \sigma^{g^0} G(x,t), \; \;\text{for} \; x \in \Omega, \; \, t \geq 0 \; \text{and} \;\sigma \geq 1.
\end{equation}
\begin{equation}
\sigma^{g^0} G(x,t) \leq G(x,\sigma t) \leq \sigma^{g_0} G(x,t), \; \;\text{for} \; x \in \Omega, \; \, t \geq 0 \; \text{and} \;\sigma \leq 1.
\end{equation}
\item There exists a constant $C>0$ such that
\begin{equation}
G(x, a + b) \leq C(G(x,a) + G(x,b))
\end{equation}
\end{enumerate}
\end{lm}

\begin{df} 
We define $G^*(\cdot)$ the conjugate $\Phi$-function of $G(\cdot)$, by
$$ G^*(x,s) := \sup_{t \geq 0}(st - G(x,t)), \; \, \; \text{for} \; x \in \Omega \; \text{and} \; s \geq 0.$$ 
Note that $G^*(\cdot)$ is also a generalized $\Phi$-function and can be represented as
$$ G^*(x,t) =  \displaystyle \int_{0}^t g^{-1}(x,s) \, \mathrm{d}s,$$
with $g^{-1}(x,s) : = \sup \{ t \geq 0 \; : \; g(x,t) \leq s \}$. 
\end{df}

\begin{lm} \label{Prs of G^*}
 Let $G(\cdot)$ satisfies $(SC)$.
\begin{enumerate}[label=(\alph*)]
\item Then $G^*(\cdot)$ satisfies also $(SC)$, as follows
$$\displaystyle \frac{g^0+1}{g^0} \leq \displaystyle \frac{t g^{-1}(x,t)}{G^*(x,t)} \leq \frac{g_0+1}{g_0}.$$
\item The functions $G(\cdot)$ and $G^*(\cdot)$ satisfies the following Young inequality 
$$ st \leq G(x,t) + G^*(x,s), \; \, \text{for} \; x \in \Omega \; \text{and} \; s,t \geq 0.$$
Further, we have the equality if $s = g(x,t)$ or  $t = g^{-1}(x,s)$.
\item The functions $G(\cdot)$ and $G^*(\cdot)$ satisfy the H\"older inequality
$$ \left| \displaystyle \int_{\Omega} u(x)v(x) \, \mathrm{d}x \right| \leq 2||u||_{G(\cdot)} ||v||_{G^*(\cdot)}, \; \; \text{for} \; u \in L^{G(\cdot)}(\Omega) \; \text{and} \; v \in L^{G^*(\cdot)}(\Omega).$$
\end{enumerate}
\end{lm}

\noindent \textbf{Musielak-Orlicz-Sobolev spaces.} The generalized Orlicz space, also called Musielak-Orlicz space,  is defined as the set
$$ L^{G(\cdot)}(\Omega) = \{ u \in L^0(\Omega) \; : \; \displaystyle \int_{\Omega} G(x,|u|) \, \mathrm{d}x < \infty \},$$
equipped with the following norms:
\begin{itemize}
\item Luxembourg norm: 
$$ \norme{u}_{G(\cdot)} := \inf \{ \lambda > 0 \; : \; \displaystyle \int_{\Omega} G\left(x,\frac{|u|}{\lambda}\right) \, \mathrm{d}x \leq 1 \}.$$
\item Orlicz norm: 
$$\norme{u}_{G(\cdot)}^0 = \sup \{ |\displaystyle \int_{\Omega} u(x)v(x) \, \mathrm{d}x| \; : \; v \in L^{G^*(\cdot)}(\Omega), \; \displaystyle \int_{\Omega} G^*\left(x,|v| \right) \, \mathrm{d}x \leq 1  \}.$$
\end{itemize}
These norms are equivalent, precisely we have
$$\norme{u}_{G(\cdot)} \leq \norme{u}_{G(\cdot)}^0 \leq 2 \norme{u}_{G(\cdot)}.$$
Then, by definition of Orlicz norm and Young inequality, we have
\begin{equation}
\norme{u}_{G(\cdot)} \leq \norme{u}_{G(\cdot)}^0 \leq \displaystyle \int_{\Omega} G(x,|u|) \, \mathrm{d}x + 1.
\end{equation}
Next, we define the Musielak-Orlicz-Sobolev space by
$$ W^{1,G(\cdot)}(\Omega) := \{ u \in L^{G(\cdot)}(\Omega) \; : \; |\nabla u| \in L^{G(\cdot)}(\Omega), \; \, \text{in the distribution sense} \},$$
equipped with the norm
$$ ||u||_{1,G(\cdot)} = ||u||_{G(\cdot)} + ||\nabla u||_{G(\cdot)}.$$

\noindent The following lemmas establish properties of convergent sequences in generalized Orlicz spaces.

\begin{lm} \label{Cvg of G}
Let $G(\cdot) \in \Phi(\Omega)$ satisfies $(SC)$, then the following relations hold true
\begin{enumerate} [label=(\alph*)]
\item $\norme{u}_{G(\cdot)}^{g_0} \leq \displaystyle \int_{\Omega} G(x,|u|) \, \mathrm{d}x  \leq \norme{u}_{G(\cdot)}^{g^0}, \; \forall u \in L^{G(\cdot)}(\Omega) \; \text{with} \; \norme{u}_{G(\cdot)} > 1.$
\item $\norme{u}_{G(\cdot)}^{g^0} \leq \displaystyle \int_{\Omega} G(x,|u|) \, \mathrm{d}x  \leq \norme{u}_{G(\cdot)}^{g_0}, \; \forall u \in L^{G(\cdot)}(\Omega) \; \text{with} \; \norme{u}_{G(\cdot)} < 1.$
\item For any sequence $\{u_i\}_i$ in $L^{G(\cdot)}(\Omega)$, we have
$$||u_i||_{G(\cdot)} \rightarrow 0  \; \; (\text{resp.} 1 ;\infty) \Longleftrightarrow \displaystyle \int_{\Omega} G(x,|u_i(x)|) \, \mathrm{d}x \rightarrow 0  \; \; (\text{resp.} 1 ;\infty).$$
\end{enumerate}
\end{lm}

\noindent In the following lemma, we have the continuous embedding into the Musielak-Orlicz spaces.

\begin{lm} \label{Em of G}
Let $G(\cdot), \psi(\cdot) \in \Phi(\Omega)$. Then $L^{G(\cdot)}(\Omega) \hookrightarrow L^{\psi(\cdot)}(\Omega)$ if and only if there exist $C>0$ and $h \in L^1(\Omega)$ with $||h||_1 < 1$ such that 
$$\psi(x,t) \leq C ( G(x,t) + h(x)),$$
for all $x \in \Omega$ and all $t \geq 0$.
\end{lm}

\noindent To study boundary value problems, we need the concept of weak boundary value spaces. We define $W^{1,G(\cdot)}_0(\Omega)$ as the closure of $C^\infty_0(\Omega)$ in $W^{1,G(\cdot)}(\Omega)$. Next, we recall the norm version of the Poincaré inequality, which will be used in this work. 

\begin{lm}\label{PC of G}
Let $\Omega$ be a bounded set of $\R^N$ and $G(\cdot) \in \Phi(\Omega)$ satisfy $(SC)$, $(A_0)$ and $(A_{1})$. For every $u \in W^{1,G(\cdot)}_0(\Omega)$, we have
$$||u||_{G(\cdot)} \leq  C||\nabla u||_{G(\cdot)}. $$
In particular, $||\nabla u||_{G(\cdot)}$ is a norm on $W^{1,G(\cdot)}_0(\Omega)$ and it is equivalent to the norm $ ||u||_{1,G(\cdot)}$.
\end{lm}

\noindent The following compact embedding theorem for Musielak-Sobolev spaces is given by Harjulehto and H\"ast\"o $\cite{ref6}$.

\begin{lm} \label{EmC of G}
Let $G(\cdot) \in \Phi(\R^N)$ satisfy $(SC)$, $(A_0)$ and $(A_{1})$ for $g^0 < N$. Suppose that $\psi \in \Phi(\R^N)$ satisfies $t^{-\alpha} G^{-1}(x,t) \approx \psi^{-1}(x,t)$, for some $\alpha \in [0,\frac{1}{N})$. Then  
$$W^{1,G(\cdot)}_0(\Omega) \hookrightarrow \hookrightarrow L^{\psi(\cdot)}(\Omega).$$
\end{lm}

\noindent \textbf{Mountain Pass Theorem in Banach spaces.}  We recall here a version of the mountain pass theorem, which was discussed by $\cite{ref1,ref12,ref10}$. We shall apply this theorem to establish critical points for finding solutions to problem $(\ref{G()-laplacian})$.  

\begin{df}
Let $X$ be a Banach spaces and let $J \in C^1(X,\R)$. We say that $J$ satisfies Palais-Smale condition $(PS)$ in $X$ if any sequence $\{u_i\}$ in $X$ such that
\begin{enumerate}
\item[(i)] $\{J(u_i)\}$ is bounded,
\item[(ii)] $J^\prime(u_i) \to 0$ as $i \to \infty$,
\end{enumerate}
has a convergence subsequence.
\end{df}

\begin{Th} \label{MP Th}
Let $X$ be a Banach spaces and let $J \in C^1(X,\R)$ satisfy the Palais-Smale condition. Assume that $J(0) = 0$, and,
\begin{enumerate}
\item[$(MP)_1$]: There exist two positive real numbers $\eta$ and $r$ such that $J(u) \geq r$ with $||u||= \eta$,
\item[$(MP)_2$]: There exists $u_1 \in X$ such that $||u_1|| > \rho$ and $J(u) < 0$.\\
\end{enumerate}
Put $$A = \{ f \in C([0,1],X) \; : \; f(0) = 0, \; f(1) = u_1\}. $$
Set $$\beta = \inf\{\max J(f([0,1])) \; : \; f \in A \}.$$
Then $\beta \geq r$ and $\beta$ is a critical value of $J$.
\end{Th}

\section{The main result and proof}

Throughout this section, let $G(\cdot) \in \Phi(\Omega)$ satisfies $(SC)$ with $g^0 < \min\{g_0^* , N \}$ and, $(A_0)$, $(A_1)$. Let $f: \Omega \times \R \to \R$ a Carathéodory function. We consider the following problem
\begin{equation} \label{G(.)-laplacian}
\left \{
   \begin{array}{rlll}
    \displaystyle -\text{div}\left(\frac{g(x,|\nabla u|)}{|\nabla u|}\nabla u\right) & = & f(x,u)  & \text{in} \; \Omega\\
    u & = & 0 & \text{on} \; \partial\Omega,
   \end{array}
   \right.
\end{equation}
In this article, we assume that $f$ satisfy the following conditions:
\begin{enumerate}
\item[•] $(f_\alpha)$: The G($\cdot$)-subcritical growth condition
$$|f(x,t)| \leq C(1+\psi(x,t))$$
with $\psi$ is density of the generalized $\Phi$-function $\Psi(\cdot)$ verifies $\Psi^{-1}(x,t) = t^{-\alpha} G^{-1}(x,t)$ with $\alpha \in (\frac{1}{g_0} - \frac{1}{g^0} , \frac{1}{N})$.
\item[•] $(f_0)$: The G($\cdot$)-superlinear at $0$ condition
$$f(x,t) = o(t^{g^0-1}), \, t \to 0 \; \text{for} \; x \in \Omega.$$
\item[•] $(AR)$: The G($\cdot$)-Ambrosetti-Rabinowitz type condition: There exist $t_0 > 0$ and $\theta > g^0$ such that
$$0 < \theta F(x,t) \leq tf(x,t) \; \text{for} \; |t| \geq t_0,$$
with $F(x,t) := \displaystyle \int_0^t f(x,s) \, \mathrm{d}s$.
\end{enumerate}

\noindent \textbf{Weak solutions.} We say that a function $u \in W_0^{1,G(\cdot)}(\Omega)$ is a weak solution of the equation $(\ref{G(.)-laplacian})$ in $\Omega$ if
$$ \int_{\Omega} \frac{g(x,|\nabla u|)}{|\nabla u|}\nabla u \cdot \nabla \varphi \, \mathrm{d}x = \int_{\Omega} f(x,u)\varphi \, \mathrm{d}x,$$
whenever  $\varphi \in W^{1,G(\cdot)}_{0}(\Omega).$\\

\noindent Note that, if $f(x,t)=0$ the existence of weak solutions to problem $(\ref{G(.)-laplacian})$ have been proved in $\cite{ref2,ref7}$. For the general case, we denote by $X$ the Musielak-Orlicz-Sobolev space $W_0^{1,G(\cdot)}(\Omega)$ and we define the functional energy corresponding to problem $(\ref{G(.)-laplacian})$ by $J : X \to \R$
$$J(u) := \int_\Omega G(x,|\nabla u|) \, \mathrm{d}x - \int_\Omega F(x,u) \, \mathrm{d}x.$$
It's well known that standard arguments imply that $J \in C^1(X,\R)$ and this derivative is
$$\langle J^\prime(u) , v \rangle = \int_\Omega \frac{g(x,|\nabla u|)}{|\nabla u|}\nabla u \cdot \nabla v  \, \mathrm{d}x  - \int_\Omega f(x,u)v  \, \mathrm{d}x.$$ 

\noindent To prove the (PS) condition in our situation, a lack of homogeneity is a major source of difficulties. So, we developed a method inspired by Lieberman's pioneering article $\cite{ref10}$, which will allow us to overcome this problem.

\begin{lm}
Suppose $(f_\alpha)$ and $(AR)$ hold, then $J$ satisfies condition $(PS)$.
\end{lm}

\noindent \textbf{Proof.} Let us assume that there exits a sequence $\{u_i\} \subset X$ such that 
$$|J(u_i)| \leq C  \quad J^\prime(u_i) \to 0.$$
Then, by the condition $(AR)$, we have
\begin{equation}
\begin{array}{ll}
C & \geq J(u_i)\\
& = \displaystyle \int_\Omega G(x,|\nabla u_i|) \, \mathrm{d}x - \int_\Omega F(x,u_i) \, \mathrm{d}x\\
& = \displaystyle \int_\Omega G(x,|\nabla u_i|) \, \mathrm{d}x - \int_{\{|u_i| < t_0\}} F(x,u_i) \, \mathrm{d}x - \int_{\{|u_i| \geq t_0\}} F(x,u_i) \, \mathrm{d}x.\\
\end{array}
\end{equation}
We will estimate separately the last two integrals. By the condition $(f_\alpha)$, we have
$$\begin{array}{ll}
\displaystyle \int_{\{|u_i| < t_0\}} F(x,u_i) \, \mathrm{d}x & \leq C \bigg(\displaystyle \int_{\{|u_i| < t_0\}} |u_i| \, \mathrm{d}x + \int_{\{|u_i| < t_0\}} \Psi(x, |u_i|) \, \mathrm{d}x \bigg)\\
& \leq C\bigg(t_0 + \displaystyle \int_{\{|u_i| < t_0\}} \Psi(x, t_0) \, \mathrm{d}x\bigg).
\end{array}$$
Or, from the structural condition of $\Psi$ (see introduction)
$$\frac{1}{-\alpha + \frac{1}{g_0}} \leq \frac{t \psi(x,t)}{\Psi(x,t)} \leq \frac{1}{-\alpha + \frac{1}{g^0}},$$
and Lemma $\ref{Prs of G}$, we have
$$ \Psi(x, t_0) \leq \max\{t_0^{-\alpha + \frac{1}{g_0}},t_0^{-\alpha + \frac{1}{g^0}}\}\Psi(x,1).$$
As $G(\cdot)$ satisfies the condition $(A_0)$, then $G^{-1}(\cdot)$ and $\psi(\cdot)$ also satisfies $(A_0)$. Hence
\begin{equation}
\displaystyle \int_{\{|u_i| < t_0\}} F(x,u_i) \, \mathrm{d}x \leq C.
\end{equation}
For the second integral, using the conditions $(AR)$ and $(f_\alpha)$, we get 
$$\begin{array}{ll}
\displaystyle \int_{\{|u_i| \geq t_0\}} F(x,u_i) \, \mathrm{d}x & \leq \displaystyle \frac{1}{\theta}  \int_{\{|u_i| \geq t_0\}} f(x,u_i)u_i \, \mathrm{d}x\\
& \leq \displaystyle \frac{1}{\theta} \bigg(\int_{\Omega} f(x,u_i)u_i \, \mathrm{d}x + \displaystyle \int_{\{|u_i| < t_0\}} |f(x,u_i)||u_i| \, \mathrm{d}x \bigg)\\
& \leq \displaystyle \frac{1}{\theta} \bigg(\int_{\Omega} f(x,u_i)u_i \, \mathrm{d}x + \displaystyle C\int_{\{|u_i| < t_0\}} (1 + \psi(x,t_0) t_0 \, \mathrm{d}x\bigg).\\
\end{array}$$
By the same method of inequality $(3.3)$, we obtain
$$\displaystyle \int_{\{|u_i| < t_0\}} (1 + \psi(x,t_0) t_0 \, \mathrm{d}x \leq C.$$
Hence,
\begin{equation}
\displaystyle \int_{\{|u_i| \geq t_0\}} F(x,u_i) \, \mathrm{d}x \leq \displaystyle \frac{1}{\theta} \int_{\Omega} f(x,u_i)u_i \, \mathrm{d}x +  C.
\end{equation}
By combining inequalities $(3.2)$, $(3.3)$ and $(3.4)$, we get
\begin{equation}
C \geq \displaystyle \int_\Omega G(x,|\nabla u_i|) \, \mathrm{d}x -  \frac{1}{\theta} \int_{\Omega} f(x,u_i)u_i \, \mathrm{d}x.
\end{equation}
In other hand, we have $ J^\prime(u_i) \to 0$ then, we can choose $i_0 \in \N$ such that for all $i \geq i_0$, we have
$$|\langle J^\prime(u_i) , v \rangle| \leq ||v||_{1,G(\cdot)} $$
So, for $v = u_i$ and the condition $(SC)$,  we obtain
$$\begin{array}{ll}
-||u_i||_{1,G(\cdot)} & \leq  \displaystyle \int_\Omega \frac{g(x,|\nabla u_i|)}{|\nabla u_i|}\nabla u_i \cdot \nabla u_i  \, \mathrm{d}x  - \int_\Omega f(x,u_i)u_i  \, \mathrm{d}x\\
& \leq g^0 \displaystyle \int_\Omega G(x,|\nabla u_i|) \, \mathrm{d}x - \int_{\Omega} f(x,u_i)u_i \, \mathrm{d}x. 
\end{array}$$
Therefore, by the inequality $(3.5)$ and Lemma $\ref{Cvg of G}$, we have
$$\begin{array}{ll}
C & \geq \displaystyle \bigg(1 - \frac{g^0}{\theta}\bigg)\displaystyle \int_\Omega G(x,|\nabla u_i|) \, \mathrm{d}x - \frac{1}{\theta}||u_i||  
\\
& \geq \displaystyle \bigg(1 - \frac{g^0}{\theta}\bigg) \min\{||u_i||_{1,G(\cdot)}^{g^0},||u_i||_{1,G(\cdot)}^{g_0} \} - \frac{1}{\theta}||u_i||.
\end{array}$$
As $\theta > g^0$, then $\{u_i\}$ is bounded in $X$. Which implies $u_i \rightharpoonup u$ in $X$. As $\psi^{-1}(x,t) = t^{-\alpha} G^{-1}(x,t)$ with $\alpha \in (\frac{1}{g_0} - \frac{1}{g^0} , \frac{1}{N})$, then $X \hookrightarrow  \hookrightarrow  L^{\psi(\cdot)}(\Omega)$, so $u_i \to u$ in $L^{\Psi(\cdot)}(\Omega)$. Thus, by H\"older inequality in $L^{\Psi(\cdot)}(\Omega)$, we have
$$\begin{array}{ll}
\displaystyle \int_{\Omega} f(x,u_i)(u_i - u) \, \mathrm{d}x & \leq \displaystyle C\int_{\Omega} (1 + \psi(x,|u_i|)|u_i - u| \, \mathrm{d}x \\
& \leq C||1 + \psi(x,|u_i|)||_{\Psi^*(\cdot)}||u_i - u||_{\Psi(\cdot)}.
\end{array}$$
So, by the inequality $(2.5)$, the Young equality, the condition $(SC)$ and Lemma $\ref{Cvg of G}$, we have
$$\begin{array}{ll}
||1 + \psi(x,|u_i|)||_{\Psi^*(\cdot)} & \leq 1 + ||\psi(x,|u_i|)||_{\Psi^*(\cdot)} \\
&\leq \displaystyle \int_{\Omega} \Psi^*(x,\psi(x,|u_i|)) \, \mathrm{d}x + 2\\
& \leq C\displaystyle \int_{\Omega} \Psi(x,|u_i|) \, \mathrm{d}x + 2\\
& \leq C \max\{||u_i||_{\Psi(\cdot)}^{\frac{1}{-\alpha + \frac{1}{g_0}}} , ||u_i||_{\Psi(\cdot)}^{\frac{1}{-\alpha + \frac{1}{g^0}}} \} + 2\\
& \leq C.
\end{array}$$
Then
$$\displaystyle \int_{\Omega} f(x,u_i)(u_i - u) \, \mathrm{d}x \leq C ||u_i - u||_{\Psi(\cdot)} \to 0.$$
As $|\langle J^\prime(u_i) , u_i - u \rangle| \to 0$, then
$$ \displaystyle \int_\Omega \frac{g(x,|\nabla u_i|)}{|\nabla u_i|}\nabla u_i \cdot (\nabla u - \nabla u_i)  \, \mathrm{d}x \to 0.$$ 
Furthermore, we have
$$ \displaystyle \int_\Omega \frac{g(x,|\nabla u|)}{|\nabla u|}\nabla u \cdot (\nabla u - \nabla u_i)  \, \mathrm{d}x \to 0.$$ 
Hence
$$ \displaystyle \int_\Omega \bigg(  \frac{g(x,|\nabla u|)}{|\nabla u|}\nabla u - \frac{g(x,|\nabla u_i|)}{|\nabla u_i|}\nabla u_i\bigg) \cdot (\nabla u - \nabla u_i ) \, \mathrm{d}x \to 0$$ 
Using the condition $(SC)$ and Cauchy-Schwarz inequality, we have for $\theta_t = tu + (1-t)u_i$, with $t \in (0,1)$.
$$\begin{array}{ll}
& \displaystyle \bigg( \frac{g(x,|\nabla u|)}{|\nabla u|}\nabla u  - \frac{g(x,|\nabla u_i|)}{|\nabla u_i|}\nabla u_i\bigg) \cdot (\nabla u - \nabla u_i)\\
& = \displaystyle \left(\int_0^1 \frac{\partial}{\partial t}\left(\frac{g(x,|\nabla \theta_t|)}{|\nabla \theta_t|}\nabla \theta_t\right)  \, \mathrm{d}t\right)\cdot (\nabla u - \nabla u_i)\\
& = \displaystyle |\nabla u - \nabla u_i|^2 \int_0^1 \frac{g(x,|\nabla \theta_t|)}{|\nabla \theta_t|}  \, \mathrm{d}t\\
& \quad + \displaystyle  \int_0^1 g(x,|\nabla \theta_t|)\left( \frac{|\nabla \theta_t| g^\prime(x,|\nabla \theta_t|)}{g(x,|\nabla \theta_t|)} - 1\right) \frac{(\nabla \theta_t \cdot (\nabla u - \nabla u_i))^2}{|\nabla \theta_t|^3}  \, \mathrm{d}t \\
& \geq  \displaystyle \min(1,g_0-1) |\nabla u - \nabla u_i|^2 \int_0^1 \frac{g(x,|\nabla \theta_t|)}{|\nabla \theta_t|}  \, \mathrm{d}t.
\end{array}$$
Which implies
$$\begin{array}{ll}
& \displaystyle \int_\Omega \bigg(  \frac{g(x,|\nabla u|)}{|\nabla u|}\nabla u - \frac{g(x,|\nabla u_i|)}{|\nabla u_i|}\nabla u_i\bigg) \cdot (\nabla u - \nabla u_i ) \, \mathrm{d}x\\
& \geq   \displaystyle \min(1,g_0) \int_{\Omega}\int_0^1 \frac{g(x,|\nabla \theta_t|)}{|\nabla \theta_t|}|\nabla u - \nabla u_i|^2  \,  \mathrm{d}t \mathrm{d}x. 
\end{array}$$
Now we write $S_1 = \{ x \in \Omega \; , \; \, |\nabla u - \nabla u_i| \leq 2|\nabla u|\}$ and $S_2 = \{ x \in \Omega \; , \; \, |\nabla u - \nabla u_i| > 2|\nabla u|\}$. Then $S_1 \cup S_2 = \Omega$ and 
$$\begin{array}{ll}
\displaystyle\frac{1}{2}|\nabla u| \leq |\nabla \theta_t| \leq 3 |\nabla u| \; \; & \text{on} \; S_1  \; \; \text{for} \; t \geq \displaystyle \frac{3}{4},\\
\\
\displaystyle \frac{1}{4}|\nabla u - \nabla u_i| \leq |\nabla \theta_t| \leq 3 |\nabla u - \nabla u_i| \; \; & \text{on}  \; S_2  \; \; \text{for} \; t \leq \displaystyle \frac{1}{4}.\\
\end{array}$$
Therefore
\begin{equation}
\begin{array}{ll}
& \displaystyle \int_\Omega \bigg(  \frac{g(x,|\nabla u|)}{|\nabla u|}\nabla u - \frac{g(x,|\nabla u_i|)}{|\nabla u_i|}\nabla u_i\bigg) \cdot (\nabla u - \nabla u_i ) \, \mathrm{d}x \\ 
& \geq C \bigg( \displaystyle \int_{S_1} \frac{g(x,|\nabla u|)}{|\nabla u|}|\nabla u - \nabla u_i|^2  \, \mathrm{d}x  + \displaystyle \int_{S_2} G(x,|\nabla u - \nabla u_i|)  \, \mathrm{d}x \bigg).
\end{array}
\end{equation}
Hence 
\begin{equation}
\displaystyle \int_{S_2} G(x,|\nabla u - \nabla u_i|)  \, \mathrm{d}x \leq  C\displaystyle \int_\Omega \bigg(  \frac{g(x,|\nabla u|)}{|\nabla u|}\nabla u - \frac{g(x,|\nabla u_i|)}{|\nabla u_i|}\nabla u_i\bigg) \cdot (\nabla u - \nabla u_i ) \, \mathrm{d}x.
\end{equation}
To estimate the integrals over $S_1$, using the condition $(SC)$, $t \to g(x,t)$ is a nondecreasing function, and the H\"older inequality in $L^2(S_1)$, we have
$$\begin{array}{ll}
& \displaystyle \int_{S_1} G(x,|\nabla u - \nabla u_i|)  \, \mathrm{d}x\\
& \leq C \displaystyle \int_{S_1} g(x,|\nabla u - \nabla u_i|)|\nabla u - \nabla u_i|  \, \mathrm{d}x\\
& \leq  C \displaystyle \int_{S_1} \left(\frac{g(x,|\nabla u|)}{|\nabla u|}|\nabla u|\right)^\frac{1}{2}|\nabla u - \nabla u_i|\left(g(x,|\nabla u - \nabla u_i|)\right)^\frac{1}{2}|  \, \mathrm{d}x\\
& \leq C \left(\displaystyle \int_{S_1} \frac{g(x,|\nabla u|)}{|\nabla u|}|\nabla u - \nabla u_i|^2  \, \mathrm{d}x\right)^\frac{1}{2}\left(\displaystyle \int_{S_1} g(x,|\nabla u|)|\nabla u|  \, \mathrm{d}x\right)^\frac{1}{2}\\
& \leq C \left(\displaystyle \int_{S_1} \frac{g(x,|\nabla u|)}{|\nabla u|}|\nabla u - \nabla u_i|^2  \, \mathrm{d}x\right)^\frac{1}{2}\left(\displaystyle \int_{S_1} G(x,|\nabla u|)  \, \mathrm{d}x\right)^\frac{1}{2}.
\end{array}$$
Hence, from the inequality $(3.6)$, we have
\begin{equation}
\begin{array}{ll}
& \displaystyle \int_{S_1} G(|\nabla u - \nabla u_i|)  \, \mathrm{d}x \\
& \leq  C \left(\displaystyle \int_\Omega \bigg(  \frac{g(x,|\nabla u|)}{|\nabla u|}\nabla u - \frac{g(x,|\nabla u_i|)}{|\nabla u_i|}\nabla u_i\bigg) \cdot (\nabla u - \nabla u_i ) \, \mathrm{d}x \right)^{\frac{1}{2}} \left(\displaystyle \int_{\Omega} G(x,|\nabla u|)  \, \mathrm{d}x\right)^\frac{1}{2}.
\end{array}
\end{equation}
Collecting the inequalities $(3.7)$ and $(3.8)$, we have
$$\begin{array}{ll}
& \displaystyle \int_{\Omega} G(|\nabla u - \nabla u_i|)  \, \mathrm{d}x\\
& \quad \leq  \displaystyle C \left( \int_\Omega \bigg(  \frac{g(x,|\nabla u|)}{|\nabla u|}\nabla u - \frac{g(x,|\nabla u_i|)}{|\nabla u_i|}\nabla u_i\bigg) \cdot (\nabla u - \nabla u_i ) \, \mathrm{d}x \right)^{\frac{1}{2}}\left(\int_{\Omega} G(|\nabla u|) \, \mathrm{d}x \right)^\frac{1}{2}\\
& \quad \quad + \displaystyle \int_\Omega \bigg(  \frac{g(x,|\nabla u|)}{|\nabla u|}\nabla u - \frac{g(x,|\nabla u_i|)}{|\nabla u_i|}\nabla u_i\bigg) \cdot (\nabla u - \nabla u_i ) \, \mathrm{d}x  \bigg) \longrightarrow 0.
\end{array}$$
Therefore, by Lemma $(\ref{Cvg of G})$ and Poincaré inequality, we have $u_i \longrightarrow u \; \; \text{in} \; X$. \qed

\begin{Th}
Assume that $f$ satisfy $(AR)$, $(f_0)$ and $(f_\alpha)$ with $\alpha \in (\frac{1}{g_0} - \frac{1}{g^0} , \frac{1}{N})$. Then problem $(\ref{G(.)-laplacian})$ has nontrivial weak solution.
\end{Th}

\noindent \textbf{Proof.} We show that $J$ satisfies all assumptions of Theorem $\ref{MP Th}$. We start with condition $(MP)_1$. For this, we have
$$J(u) = \int_\Omega G(x,|\nabla u|) \, \mathrm{d}x - \int_\Omega F(x,u) \, \mathrm{d}x.$$
Let $||u||_{1, G(\cdot)} < 1$. Then
$$J(u) \geq ||\nabla u||_{G(\cdot)}^{g^0} - \int_\Omega F(x,u) \, \mathrm{d}x. $$
From conditions $(f_0)$, $(f_\alpha)$ and condition $(A_0)$, we have
$$\begin{array}{ll}
F(x,t) & \leq \epsilon |t|^{g^0}\mathbbm{1}_{\{|s| \leq 1\}}(t) + C (t + \Psi(x,t))\mathbbm{1}_{\{|s| > 1\}}(t)\\
 & \leq \epsilon |t|^{g^0}\mathbbm{1}_{\{|s| \leq 1\}}(t) + C \Psi(x,t)\mathbbm{1}_{\{|s| > 1\}}(t).
\end{array}$$
Then 
$$\int_\Omega F(x,u) \, \mathrm{d}x \leq \epsilon \int_{\{|u| \leq 1 \}} |u|^{g^0} \, \mathrm{d}x + C \int_{\{|u| > 1 \}}\Psi(x,|u|) \, \mathrm{d}x.$$
By the condition $(A_0)$, we have $t^{g^0} \leq C G(x,t)$ for $t \leq 1$. Then, using Lemma $\ref{Em of G}$ , we get 
$$ \int_{\{|u| < 1 \}} |u|^{g^0} \, \mathrm{d}x \leq ||u||^{g^0}_{g^0} \leq C ||u||^{g^0}_{1, G(\cdot)} \leq C ||\nabla u||^{g^0}_{G(\cdot)}$$
Furthermore, by Lemma $\ref{Cvg of G}$ and compact embedding in the Musielak-Orlicz-Sobolev spaces, we have
$$ \int_{\Omega} \Psi(x,|u|) \, \mathrm{d}x \leq C ||u||_{\Psi(\cdot)}^{\frac{1}{-\alpha + \frac{1}{g_0}}} \leq C ||\nabla u||_{G(\cdot)}^{\frac{1}{-\alpha + \frac{1}{g_0}}}.$$
Therefore
$$J(u) \geq ||\nabla u||_{G(\cdot)}^{g^0} - \epsilon C||\nabla u||_{G(\cdot)}^{g^0} - ||\nabla u||_{G(\cdot)}^{\frac{1}{-\alpha + \frac{1}{g_0}}}$$
So, choosing $\epsilon = \frac{1}{2C}$, we obtain
$$J(u) \geq \frac{1}{2}||\nabla u||_{G(\cdot)}^{g^0} - ||\nabla u||_{G(\cdot)}^{\frac{1}{-\alpha + \frac{1}{g_0}}}$$
As $\alpha > \displaystyle  \frac{1}{g_0} - \frac{1}{g^0}$ which implies that $g^0 < \displaystyle  \frac{1}{-\alpha + \frac{1}{g_0}}$, then there exist two constants $\eta , r$ such that $J(u) \geq r > 0$ with $||\nabla u||_{G(\cdot)} = \eta \in (0 , 1)$.\\
For the second condition $(MP)_2$, We use an important consequence of the condition $(AR)$:
$$F(x,t) \geq C|t|^\theta \; \text{for} \; |t| \geq t_0. $$
So, for $\omega \in X$, $\omega>0$ and $t > 1$, let us denote
$$M_t(\omega) = \{x \in \Omega \; : \; t\omega(x) \geq t_0 \} $$
We choose $u \in X$, $u >0$ with $|M_t(u)| > 0$. It is clear that $M_1(u) \subset M_t(u)$ and hence $|M_1(u)| \leq |M_t(u)|$ for all $t > 1$. Then
$$\int_{M_t(u)} F(x,tu) \, \mathrm{d}x \geq t^\theta C \int_{M_1(u)} u^\theta \, \mathrm{d}x. $$
Furthermore, by the condition $(f_\alpha)$ and the condition $(A_0)$, we have
$$\begin{array}{ll}
\displaystyle  \int_{M_t(u)^c} F(x,tu) \, \mathrm{d}x & \leq  C \displaystyle  \int_{M_t(u)^c} tu + \Psi(x,tu) \, \mathrm{d}x \\
& \leq  C \displaystyle  \int_{M_t(u)^c} t_0 + \Psi(x,t_0) \, \mathrm{d}x \\
& \leq  C
\end{array}$$
Hence
$$ J(tu) \geq t^{g_0} \int_\Omega G(x,|\nabla u|) \, \mathrm{d}x - t^\theta C \int_{M_1(u)} u^\theta \, \mathrm{d}x - C.$$
Since $\theta > g^0$ then $J(tu) \to - \infty$ when $t \to \infty$. The fact $J(0)=0$, $J$ satisfies the all assumptions of Theorem $\ref{MP Th}$. Therefore $J$ has at least one nontrivial critical point, i.e problem $(\ref{G(.)-laplacian})$ has a nontrivial weak solution. The proof is complete. \qed

\end{document}